\documentclass[12pt,reqno]{amsart}
\usepackage{amssymb}
\usepackage{amsmath}

\usepackage{amssymb, amsmath, amsfonts, latexsym}
\usepackage{enumerate}

\def\Lip{{\mathrm{{\rm Lip}}}}

\newtheorem{thm}{Theorem}[section]
\newtheorem{cor}[thm]{Corollary}
\newtheorem{lem}[thm]{Lemma}
\newtheorem{prop}[thm]{Proposition}
\newtheorem{example}[thm]{Example}
\newtheorem{remarks}[thm]{Remark}

\theoremstyle{definition}
\newtheorem{defn}{Definition}[section]

\numberwithin{equation}{section} \theoremstyle{remark}

\title[Log-Sobolev, isoperimetry and  transport inequalities  on graphs]{Log-Sobolev, isoperimetry and  transport inequalities  on graphs}

\author{Yutao Ma}
\address{Yutao MA\\ School of Mathematical Sciences $\&$ Lab. Math. Com. Sys., Beijing Normal University, 100875, Beijing China.
Partially supported by NSFC 11431014, 11371283, YETP0264, 985 Projects and the Fundamental Research Funds for the Central Universities.
} \email{mayt@bnu.edu.cn}

\author{Ran Wang}
\address{Ran Wang \\School of Mathematical sciences, University of Science and Technology of China, No.96, JinZhai Road Baohe District, Hefei, Anhui, 230026, P.R.China. 
Partially supported  by NSFC 11301498, 11431014 and the Fundamental Research Funds for the Central Universities WK0010000048. 
} \email{wangran@ustc.edu.cn}

\author{Liming Wu}
\address{Liming Wu\\ Institute of Applied Math., Chinese Academy of Sciences,
100190, Beijing, China and Laboratoire de Math. CNRS-UMR 6620,
Universit\'e Blaise Pascal, 63177 Aubi\`ere, France.
Partially supported by Thousand Talents Project and le projet ANR EVOL.
}
\email{Li-Ming.Wu@math.univ-bpclermont.fr}

\newcommand{\ee}{\mathbb{E}}

\newcommand{\rr}{\mathbb{R}}
\newcommand{\pp}{\mathbb{P}}

\newcommand{\zz}{\mathbb{Z}}

\def\Lip{{\mathrm{{\rm Lip}}}}

\def\FF{\mathcal F}
\def\EE{\mathcal E}

\def\LL{\mathcal L}
\def\MM{\mathcal M}

 \def\e{\mathbf{e} }

\def\<{\langle}
\def\>{\rangle}

\def\beq{\begin{equation}}
\def\nneq{\end{equation}}

\def\bdef{\begin{defn}}
\def\ndef{\end{defn}}

\def\bthm{\begin{thm}}
\def\nthm{\end{thm}}

\def\bprop{\begin{prop}}
\def\nprop{\end{prop}}

\def\brmk{\begin{remarks}}
\def\nrmk{\end{remarks}}

\def\bexa{\begin{example}}
\def\nexa{\end{example}}

\def\blem{\begin{lem}}
\def\nlem{\end{lem}}

\def\bcor{\begin{cor}}
\def\ncor{\end{cor}}

\date{}

\def\bexe{\begin{exe}}
\def\nexe{\end{exe}}

\def\bprf{\begin{proof}}
\def\nprf{\end{proof}}

\def\bdes{\begin{description}}
\def\ndes{\end{description}}

\begin{document}
\maketitle

\begin{abstract} In this paper, we study some functional inequalities   (such as Poincar\'e inequalities,  logarithmic Sobolev inequalities, generalized Cheeger isoperimetric inequalities,
transportation-information inequalities and transportation-entropy
inequalities) for reversible nearest-neighbor Markov processes on a
connected finite graph by means of  (random) path
method. We provide estimates of the involved
constants.
\end{abstract}

\medskip
\noindent {\bf MSC 2010 :}  60E15; 05C81; 39B72.

\medskip
\noindent{\bf Keywords :} graph; log-Sobolev inequality;
isoperimetry; transport inequality.

\section {Introduction}

Let $G=(V, E)$ be a finite connected graph with vertex set $V$ and
oriented edges set $E$, which is a symmetric subset of
$V^2\backslash \{(x,x);x\in V\}$. If $(x,y)\in E$, we call that
$x,y$ are adjacent, denoted by $x\sim y$. Consider the operator
\beq\label{LL} \LL f(x)=\sum_y q(x,y)\left(f(y)-f(x)\right), \ \
\text{for all} \ x \in V \nneq for any function $f : V\rightarrow
\rr$, where the jump rate $q(x,y)$ from $x$ to $y$ is non-negative, and
$q(x,y)>0$ if and only if $x\sim y$.

Let $(X_t)$ be the Markov process generated by $\LL$, defined on
$(\Omega, (\FF_t)_{t\ge0}, (\pp_x)_{x\in V})$. We always assume the
{\it reversibility} condition, i.e. there is some probability
measure $\mu$ satisfying the {\it detailed balance condition}
\beq\label{balance} Q(x,y):=\mu(x)q(x,y)=\mu(y)q(y,x), \ \forall
(x,y)\in E. \nneq Equivalently, the operator $\LL$ is self-adjoint
on $L^2(\mu)$, that is
$$\aligned
\langle f,-\LL g\rangle_{\mu}&=\langle -\LL f,
g\rangle_{\mu}=\frac{1}{2}\sum_{x,y}\left(f(x)-f(y)\right)\left(g(x)-g(y)\right)Q(x,y)\\
&=\frac 12 \sum_{\e\in E} D_{\e}f D_{\e}g Q(\e)=:\EE(f,g),
\endaligned$$
where $D_{\e}f:=f(y)-f(x)$ for $\e=(x,y)\in E$. When $q(x,y)=1/d_x$
with $d_x$ the {\it degree} of $x$ (the number of neighbors
$y\sim x$), $\LL$ becomes the {\it Laplacian $\Delta$} on the graph.
In that case $\mu(x)=d_x/|E|$ and $Q(\e)=1/|E|$.

Define the variance of $f$ with respect to $\mu$
 $$
 {\rm Var}_{\mu}(f)=\mu((f-\mu(f))^2),
 $$
and the entropy of $f^2$ under $\mu$
   $$ {\rm Ent}_\mu(f^2)=\mu(f^2\log f^2)- \mu(f^2) \log
\mu(f^2).$$

We say that $\mu$ satisfies a Poincar\'e inequality if there exists a constant $\lambda>0$ such that for all $f\in L^2(\mu)$,
\beq\label{eq Poin}
 {\rm Var}_\mu(f) \le  \lambda\EE(f,f),
\nneq
$\mu$ satisfies a  log-Sobolev
inequality if there exists a constant $\alpha>0$ such that for all $f\in L^2(\mu)$,
\beq\label{eq LS}
{\rm Ent}_\mu(f^2)\le 2\alpha\EE(f,f).
\nneq
  The optimal constants $ \lambda$ and $\alpha$ in \eqref{eq Poin} and \eqref{eq LS} are called respectively the Poincar\'e constant and the log-Sobolev constant of $\mu$, which are denoted by $c_{\rm P}$ and $c_{\rm LS}$ respectively. It is well known that $c_{\rm P}\le c_{\rm LS}$, see \cite{Rot81}.
\vskip0.3cm
 The Poincar\'e inequality and logarithmic Sobolev inequality play a crucial role in the analysis of the behaviour of the process.
  To study of the Poincar\'e constant, the path combinatoric method was introduced by Jerrum and Sinclair \cite{JS} in theoretic computer science, and further developed  by Diaconis and Stroock \cite{DS}, Fill \cite{F},   Sinclair \cite{Sin},  Chen \cite{Chen99}, and so on.
The logarithmic Sobolev inequality  in the discrete setting was also  studied by many authors, such as Diaconis
and Saloff-Coste \cite{D-SC}, Roberto \cite{Rob},  Lee and Yau \cite{LY}, Chen \cite{Chen00}, Chen and Sheu \cite{CC03}, Chen {\it et al.} \cite{CLS}, and so on.
The reader is referred to the books of Saloff-Coste \cite{SC} and Chen \cite{Chen05} for further
information.
\vskip0.3cm

The main purpose of this paper is to study  the logarithmic Sobolev
inequality, the generalized Cheeger isoperimetric inequality and the transport inequality.

\vskip0.3cm
The remainder of this paper is organized as follows: in next section, we focus on the logarithmic Sobolev inequality,
the third section is devoted to the transportation-information inequality and the generalized Cheeger isoperimetric inequality.  In the last section,
some examples are discussed and  the estimates of involved constants are given.

\section{Logarithmic Sobolev inequality}
\subsection{Length functions, random paths}

A path $\gamma_{xy}$ from $x$ to $y$ is a family of edges
$\{\e_1,\cdots, \e_n\}$ where $\e_k=(x_{k-1},x_k)\in E$, such that
$x_0=x,x_n=y$. It is said to have no circle if all $x_k,
k=0,\cdots, n,$ are different. A positive function $w:E\rightarrow
(0,+\infty)$ defined on the edge set $E$ is called {\it length
function}, if $w(x,y)=w(y,x)$ for any $\e=(x,y)\in E$. Given the length function $w$,  the {\it $w$-length} of a
path $\gamma_{xy}$ from $x$ to $y$ is defined by
$$ |\gamma_{xy}|_w:=\sum_{\e\in \gamma_{xy}}w(\e),
$$
and the distance associated with $w$ is
$$
\rho_w(x,y):=\min_{\gamma_{xy}} |\gamma_{xy}|_w.
$$
When $w\equiv 1$, $\rho_w=:\rho_1$ is the natural graph distance on
$V$.

\vskip0.3cm
Diaconis and Stroock \cite{DS} showed that the Poincar\'e constant
$c_{\rm P}$ satisfies that
\beq\label{10} c_{\rm P}\le \max_{\e\in E} \sum_{x,y\in V}
1_{\gamma_{xy}}(\e) |\gamma_{xy}|_{1/Q}\mu(x)\mu(y),
\nneq
for any collection of paths $\{\gamma_{xy}, x,y\in V\}$, where $|\gamma_{xy}|_{1/Q}$ is the length $|\gamma_{xy}|_w$ with
$w(\e)=1/Q(\e)$, a quite natural distance associated with the Markov
process. Furthermore, by using the  length functions, the estimate \eqref{10} can be improved to be
\beq\label{11} c_{\rm P}\le \max_{\e\in E} \frac1{Q(\e)w(\e)} \sum_{x,y\in V}
1_{\gamma_{xy}}(\e) |\gamma_{xy}|_{w}\mu(x)\mu(y), \nneq
 which is sharp for birth-death processes (see Kahale \cite{ka} or Chen \cite{Chen99}).

\vskip0.3cm

Now for any $x,y\in V$ different,  let $\gamma_{xy}$ be a random
(maybe deterministic) path without circle from $x$ to $y$. By
convention we set $\gamma_{xx}=\emptyset$, and denote by
$\ee^\gamma$ the expectation w.r.t. $\{\gamma_{xy}, x,y\in V\}$.

\subsection{Logarithmic Sobolev inequality}

\bthm\label{thm11} For any length function $w$ and any edge $\e$, let
\beq\label{thm11a} L_{w,\e}(x):=\ee^\gamma \sum_{y\in V}
1_{\gamma_{xy}}(\e)|\gamma_{xy}|_w\mu(y).\nneq
 The logarithmic Sobolev constant $c_{\rm LS}$ is bounded by
\beq\label{thm11b}c_{\rm LS} \le \inf_w
\max_{\e\in E} \frac 1{Q(\e)w(\e)}\left({\rm Ent}_\mu(L_{w,\e})+
\mu(L_{w,\e})\log(e^2+1)\right), \nneq
 where $\inf_w$ is taken over all length functions $w$ on $E$ and $e$ is the Euler constant.
  \nthm

The upper bound in (\ref{thm11b}) gives us a very practical
criterion for the logarithmic Sobolev inequality,
following the classical idea of Lyapunov function method for
stability.

The logarithmic Sobolev inequality above is based on the following weighted
Poincar\'e inequality, which is a slight generalization of \eqref{11}.

 \blem[{\bf Weighted Poincar\'e
inequality}]\label{lem12} Let $\varphi$ be a nonnegative function on
$V$, then for any length function $w$,
\beq\label{lem12a}
 \sum_{x\in V} (f(x)-\mu(f))^2 \varphi(x)\mu(x)\le c(\varphi,w)\EE(f,f), \; \forall
f:V\rightarrow \rr,
\nneq
where $$\aligned c(\varphi,w)&:=\max_{\e\in E}
\frac 2{Q(\e)w(\e)} \ee^\gamma \sum_{x,y\in
V} 1_{\gamma_{xy}}(\e)|\gamma_{xy}|_w\varphi(x)\mu(x)\mu(y)\\
&=\max_{\e\in E} \frac 2{Q(\e)w(\e)} \sum_{x\in V}
L_{w,\e}(x)\varphi(x)\mu(x).
\endaligned$$ \nlem

When $\varphi\equiv 1$, our constant $c(\varphi,w)$ is twice of
the quantity at the right hand side  of  \eqref{11}.

\bprf For any fixed realization of random paths $\{\gamma_{xy};
x,y\in V\}$,
$$
\aligned
&\sum_{x\in V} (f(x)-\mu(f))^2 \varphi(x)\mu(x)\\
=&\sum_{x\in V} \left(\sum_{y\in
V}(f(x)-f(y))\mu(y)\right)^2 \varphi(x)\mu(x)\\
=&\sum_{x\in V} \left(\sum_{y\in
V}\mu(y)\sum_{\e\in \gamma_{xy}}D_{\e}f \right)^2 \varphi(x)\mu(x)\\
\le& \sum_{x,y\in V} \left(\sum_{\e\in \gamma_{xy}}D_{\e}f \right)^2 \varphi(x)\mu(x)\mu(y)\\
\le& \sum_{x,y\in V} \left(\sum_{\e\in \gamma_{xy}} w(\e) \right)  \left(\sum_{\e\in \gamma_{xy}} \frac 1{w(\e)}(D_{\e}f)^2 \right) \varphi(x)\mu(x)\mu(y)\\
=&\frac 12\sum_{\e\in E} Q(\e) (D_{\e}f)^2\cdot \frac 2{Q(\e)w(\e)}
\sum_{x,y\in V}
1_{\gamma_{xy}}(\e)|\gamma_{xy}|_w\varphi(x)\mu(x)\mu(y),
\endaligned
$$
where the Cauchy-Schwarz inequality is applied twice.  Taking the
expectation $\ee^\gamma$ w.r.t. the randomness of $\gamma$, we get
the desired result.  \nprf

Now recall two important lemmas: the first one is due to Rothaus \cite{Rot} and the second was given by Barthe-Roberto \cite{BaRo}.

\blem\label{lem13}{\rm For any real function $f$ on $V$ and any constant
$a\in \rr$,
$$
{\rm Ent}_\mu(f^2) \le {\rm Ent}_\mu\left((f-a)^2\right) + 2
\mu\big((f-a)^2\big).
$$}
\nlem

\blem\label{lem14} {\rm For any real function $f$ on $V$,
$$
{\rm Ent}_\mu(f^2) + 2 \mu(f^2)\le \sup\left\{ \mu(f^2\varphi); \
\varphi\ge 0, \mu\left(e^\varphi\right)\le e^2+1\right\}.
$$
Consequently, by  Donsker-Varadhan's variational formula (see \cite{DV}),
$$
\mu(f^2\varphi) -\mu(f^2) \log \mu(e^\varphi)  \le {\rm Ent}_\mu(f^2), \ \ \forall \varphi,
$$
we have \beq\label{lem14b}
 \sup\left\{ \mu(f^2\varphi); \ \varphi\ge
0, \mu\left(e^\varphi\right)\le e^2+1\right\}\le {\rm Ent}_\mu(f^2)+
\mu(f^2)\log(e^2+1).\nneq}
\nlem

 \bprf[Proof of Theorem \ref{thm11}] By Lemma \ref{lem13},
$$
\aligned {\rm Ent}_\mu(f^2) &\le {\rm Ent}_\mu\left((f-\mu(f))^2\right) + 2 \mu\big((f-\mu(f))^2\big)\\
&\le \sup \left\{\sum_{x\in V} (f(x)-\mu(f))^2\varphi(x)\mu(x); \ \varphi\ge
0,
\mu\left(e^\varphi\right)\le e^2+1\right\}\\
&\le \sup_{\varphi\ge 0, \mu\left(e^\varphi\right)\le e^2+1}
c(\varphi,w)\cdot \EE(f, f),
\endaligned
$$
where the last inequality follows from \eqref{lem12a}. Moreover by (\ref{lem12a}) and
(\ref{lem14b}), we have
$$
\aligned
 \sup_{\varphi\ge 0, \mu\left(e^\varphi\right)\le e^2+1}
c(\varphi,w)&= \max_{\e\in E} \frac 2{Q(\e)w(\e)} \sup_{\varphi\ge 0,
\mu\left(e^\varphi\right)\le e^2+1}\mu(L_{w,\e} \varphi )\\
&\le \max_{\e\in E} \frac 2{Q(\e)w(\e)} \left({\rm Ent}_\mu(L_{w,\e})+
\mu(L_{w,\e})\log(e^2+1)\right),
\endaligned
$$
where (\ref{thm11b}) follows. \nprf

\section{Transportation inequalities}

In this section, we shall establish the transportation-information
inequality $W_1I$ and as a corollary, the transportation-entropy
inequality $W_1H$. For this purpose, let us introduce some notions and
notations.

\subsection{Wasserstein distance, entropy and information}
Given a metric $\rho$ on $V$, the Lipschitzian norm of a function
$g$ is denoted by $\|g\|_{\rm Lip(\rho)}$. For two probability
measures $\nu,\mu$ on $V$, say $\nu,\mu\in\MM_1(V)$,
\begin{itemize}
\item[(i)] their Wasserstein distance $W_{1,\rho}(\nu,\mu)$ associated with $\rho$  is defined as
$$
W_{1,\rho}(\nu,\mu)=\inf_{\pi}\iint_{V^2} \rho(x,y)\pi(dx,dy),
$$
where $\pi$ runs over all couplings of $(\nu,\mu)$, i.e.,
probability measures on $V^2$ such that $\pi(A\times V)=\nu(A)$ and
$\pi(V\times A)=\mu(A)$ for all Borel subsets $A$ of $V$. If
$\rho(x,y)=1_{x\ne y}$ is the discrete metric,
$W_{1,\rho}(\nu, \mu)=\frac{1}{2}\|\nu-\mu\|_{\rm TV}$ where $\|\tilde
\nu\|_{\rm TV}=\sup_{|f|\le 1}\tilde \nu(f)$ is the total variation
of a signed measure $\tilde \nu$.
\item[(ii)] The relative entropy of $\nu$ w.r.t. $\mu$ is given by
\begin{equation*}H(\nu|\mu)=\begin{cases}

         \sum_{x\in V} \nu(x) \log \dfrac{\nu(x)}{\mu(x)}, &\text{if } \nu\ll\mu; \\
               +\infty,  &\text{otherwise}.
                          \end{cases}
                          \end{equation*}
\item[(iii)] Fisher-Donsker-Varadhan information of a probability $\nu=h^2\mu$ w.r.t. $\mu$ is defined by
$$
I(\nu|\mu):=\frac{1}{2}\sum_{x,y\in V}\left(h(x)-h(y)\right)^2
Q(x,y)=\frac{1}{2}\sum_{\e\in E}(D_{\e}h)^2Q(\e),$$ where $D_{\e}h=h(y)-h(x)$
for the oriented edge $\e=(x,y)\in E$.
 \end{itemize}

\subsection{Transportation-information inequality}
Guillin {\it et al.} \cite{GLWY} introduced the following
transportation-information inequality for the given metric $\rho$,

\beq\label{WI-1} W^2_{1,\rho}(\nu,\mu)\le 2c_{\rm G} I(\nu|\mu),\
\forall \nu\in\MM_1(V), \nneq where $c_{\rm G}$ is the best constant.
In \cite{GLWY}, it is proved that \eqref{WI-1} is equivalent to
the following Gaussian concentration inequality: for all
probabilities $\nu\ll\mu$ and $\rho$-Lipschitzian function $g$ on
$V$, \beq\label{Gaussian concentration}
\pp_{\nu}\left(\frac{1}{t}\int_0^t g(X_s)ds>\mu(g)+ r\right)\leq
\left\|\frac{d\nu}{d\mu}\right\|_{L^2}\exp\left\{{-\frac{tr^2}{2c_{\rm
G}\|g\|^2_{\Lip(\rho)}}}\right\},\,\,\forall\,t,\, r
>0.
\nneq Here $(X_t)$ is the Markov process generated by $\LL$, defined
on some probability space $(\Omega, \FF,\pp)$ with initial
distribution $\nu$. So we call $c_{\rm G}$  the Gaussian concentration
constant for $(X_t)$ (associated to the metric $\rho$).

The reader is referred to the book of Villani \cite{V} for optimal transport, transport inequalities and related bibliographies.

\bthm\label{Thm WI} The transportation-information inequality
(\ref{WI-1}) holds with  \beq\label{WI} c_{\rm G}\le K:=\inf_w K(w), \nneq where the infimum is taken over all length functions $w$ and the
geometric constant $K(w)$ is given by
\beq\label{K} K(w)=\max_{\e\in
E}\frac{1}{Q(\e)w(\e)}\sum_{x,y\in V}\ee^\gamma
1_{\gamma_{xy}}(\e)\rho^2(x,y)|\gamma_{xy}|_w\mu(x)\mu(y).
\nneq
\nthm

\brmk{\rm When $\rho(x,y)=1_{x\ne y}$ (the discrete metric), $K$ coincides with the quantity in \eqref{11}. By Guillin {\it et al.} \cite[Theorem 3.1]{GLWY}, the transportation-information inequality w.r.t. the discrete metric
and the
Poincar\'e inequality are equivalent:
$$\frac {c_{\rm P}}{8}\le c_{\rm G}\le 2 c_{\rm P}.
$$
If we apply this result together with \eqref{11}, we obtain only $c_{\rm G}\le 2K$. Since $c_{\rm P}=K$  for birth-death processes (see \cite{ka,Chen99}), we get  $c_{\rm G}\ge K/8$. In other words, our estimate of $c_{\rm G}$ is of correct order. } \nrmk

  \bprf[Proof of Theorem \ref{Thm WI}]
For each probability measure $\nu=h^2\mu$, by
Kantorovich-Robinstein's identity (see \cite{V}) and the Cauchy-Schwarz
inequality,
\begin{align*}
W_{1,\rho}(\nu,\mu)=&\sup_{\|g\|_{\Lip(\rho)\le1}}\sum_{x\in V} g(x)\left(h^2(x)-1\right)\mu(x)\\
=&\frac{1}{2}\sup_{\|g\|_{\Lip(\rho)\le1}}\sum_{x,y\in V}\left(g(x)-g(y)\right)\left(h^2(x)-h^2(y)\right)\mu(x)\mu(y)\\
\le& \frac{1}{2}\sup_{\|g\|_{\Lip(\rho)\le1}}\sqrt{\sum_{x,y\in V}(g(x)-g(y))^2\left(h(x)-h(y)\right)^2\mu(x)\mu(y)}\\
 &\cdot\sqrt{\sum_{x,y}\left(h(x)+h(y)\right)^2\mu(x)\mu(y)}\\
\le &\sqrt{\sum_{x,y\in
V}\rho^2(x,y)\left(h(x)-h(y)\right)^2\mu(x)\mu(y)}.
\end{align*}
For any fixed random paths $\{\gamma_{xy}=\gamma_{xy}(\omega); x,y\in
V\}$ and the  length function $w$, we have by the Cauchy-Schwarz inequality,
\begin{align*}
&\sum_{x,y}\rho^2(x,y)\left(h(x)-h(y)\right)^2\mu(x)\mu(y)= \sum_{x,y}\rho^2(x,y)\left(\sum_{\e\in \gamma_{x,y}}D_{\e}h\right)^2\mu(x)\mu(y)\\
\le
&\sum_{x,y}\rho^2(x,y)\mu(x)\mu(y)\left(\sum_{\e\in\gamma_{xy}}(D_{\e}h)^2\frac{1}{w(\e)}\right)
\left(\sum_{\e\in\gamma_{xy}}w(\e)\right)\\
=&\sum_{\e\in E} (D_{\e}h)^2Q(\e)\cdot\frac{1}{Q(\e)w(\e)}\sum_{x,y\in V}
1_{\gamma_{x,y}}(\e)\rho^2(x,y)|\gamma_{xy}|_{w}\mu(x)\mu(y).
\end{align*}
Taking first the expectation $\ee^\gamma$ and then the maximum of the last term
over all oriented edges $\e$, we
get $c_{\rm G}\le K(w)$, the desired result. \nprf

\bcor\label{cor23} {\rm Assume that there exists  some constant $M>0$ such that
\beq\label{cor23a} \sup_{x\in V}\frac 12 \sum_{y\sim
x}\rho^2(x,y)q(x,y)\le M.\nneq Then the following transportation-entropy inequality
holds \beq\label{T1} W^2_{1,\rho}(\nu,\mu)\le \sqrt{2KM}H(\nu|\mu),\
\ \forall \nu\in \MM_1(V), \nneq or equivalently for any Lipschitzian
function $g$,
\beq\label{cor23b} \int e^{\lambda(g-\mu(g))} d\mu\le
\exp\left(\frac{\lambda^2 \sqrt{2MK}}{4}\|g\|^2_{\Lip(\rho)}\right), \
\lambda\in\rr. \nneq }\ncor

\bprf  The transportation-entropy
inequality (\ref{T1}) follows from the transportation information
inequality $(\ref{WI})$ under the condition (\ref{cor23a}), by
Guillin {\it et al.} \cite[Theorem 4.2]{GJLW}. The equivalence between (\ref{T1})
and the Gaussian concentration (\ref{cor23b}) is the famous
Bobkov-G\"{o}tze's characterization in \cite{BG}.\nprf

\bcor\label{cor24} For the Laplacian $\LL=\Delta$ on the connected
graph $G=(V,E)$, we have for the graph metric $\rho_1$,
$$
c_{\rm G}\le K\le \frac{ d_*^2 b D^3}{|E|},
$$
where $d_*=\max_{x\in V} d_x$, $D$ is the diameter of $G$ and
\beq\label{b}
b=\max_{\e}\sharp\{\rho_1\text{-shortest paths} \   \gamma: \e\in \gamma\}.
\nneq
 \ncor
\bprf Choose $\gamma_{xy}$ distributed uniformly on all shortest
paths from $x$ to $y$ and $w=1$, we see that  $K$  is bounded from
above by (noting that $Q(\e)=1/|E|, \mu(x)=d_x/|E|$)
$$
\max_{\e\in E}|E| \ee^\gamma \sum_{x,y\in V} 1_{\gamma_{xy}}(\e)
\rho_1(x,y)^3
 \left(\frac{d_*}{|E|}\right)^2\le  \frac{d_*^2 b  D^3}{ |E|}.
$$
\nprf

\subsection{Generalized Cheeger isoperimetric inequality}
Consider the following generalized
Cheeger isoperimetric inequality \beq\label{cheeger}
W_{1,\rho}(f\mu,\mu)\le\frac{c_{\rm I}}{2}\sum_{\e\in E} |D_{\e}f|Q(\e), \
\forall \nu=f\mu\in\MM_1(V),\nneq where  $c_{\rm I}$ is the best constant, called as {\it Cheeger constant} w.r.t.
the metric $\rho$.

Define the geometric constant $\kappa$
\beq\label{kappa}\kappa:=\max_{\e\in E} \frac 1{Q(\e)}\ee^\gamma
\sum_{x,y\in V} 1_{\gamma_{xy}}(\e)\rho(x,y)\mu(x)\mu(y).\nneq
\bthm\label{thm-Cheeger} It holds that
$$
c_{\rm I}\le \kappa.
$$

 \nthm \bprf By
Kantorovich-Robinstein's identity, we have for any fixed random
paths $\{\gamma_{xy}; x,y\in V\}$,
\begin{align*}
W_{1,\rho}(f\mu,\mu)&=\sup_{\|g\|_{\Lip(\rho)\le 1}}\sum_{x\in V} g(x)(f(x)-1)\mu(x)\\
&=\frac{1}{2}\sup_{\|g\|_{\Lip(\rho)\le 1}}\sum_{x,y\in V} \left(g(x)-g(y)\right)\left(f(x)-f(y)\right)\mu(x)\mu(y)\\
&\le \frac{1}{2}\sum_{x,y\in V}\rho(x,y)\sum_{\e: \e\in\gamma_{xy}}|D_{\e}f|\mu(x)\mu(y)\\
&=\frac{1}{2}\sum_{\e\in E} |D_{\e}f|Q(\e)\cdot
\frac{1}{Q(\e)}\sum_{x,y\in V}
1_{\gamma_{xy}}(\e)\rho(x,y)\mu(x)\mu(y).
\end{align*}
Taking the expectation $\ee^\gamma$, we obtain the desired result. \nprf

 \bcor[Weighted $L^1$-Poincar\'e inequality]\label{cor25} Given a positive function $\varphi$ on $V$, we have
 for any function $f$ on $V$,
\beq\label{TV} \int_V |f-\mu(f)| \varphi d\mu\le \frac
{\kappa}{2}\sum_{\e\in E} |D_{\e}f|Q(\e),
 \nneq where
$\kappa$ is given by (\ref{kappa}) with $\rho(x,y)=1_{x\ne
y}(\varphi(x)+ \varphi(y)), x,y\in V$. \ncor

\bprf Considering $(f-c_1)/c_2$ if necessary, we assume
without loss of generality that $f>0$ and $\mu(f)=1$. In that case,
for the metric $\rho(x,y)=1_{x\ne y}(\varphi(x)+ \varphi(y))$, it is
known that
$$W_{1,\rho}(f\mu,\mu)=\|\varphi(\nu-\mu)\|_{TV}=\int_V |f-1|\varphi d\mu.$$
It remains to apply Theorem \ref{thm-Cheeger}. \nprf

\brmk{\rm  Taking $\varphi\equiv 1$ and considering the
corresponding geometric constant $\kappa$,
 $(\ref{TV})$ is
equivalent to (by Bobkov and Houdr\'{e} \cite[Theorem 1.1]{BH})
 $$
 2\mu(A)\mu(A^c)\le \kappa \sum_{\e\in\partial A} Q(\e),\ \ \forall A\subset V,
  $$
where  $$\partial A:=\{\e=(x,y)\in E; x\in A, \; y\in A^c\}$$ is the
boundary of $A$. Thus $(\ref{TV})$ implies the standard Cheeger
inequality $$\mu(A)\le \kappa \sum_{\e\in\partial
A} Q(\e),\ \ \forall A\subset V \ \text{such that} \ \mu(A)\le
\frac{1}{2},$$  which has an equivalent functional version as : for every function $f$ on $V$,
\beq\label{func. Ch}
\sum_x|f(x)-{\rm med}_{\mu}(f)|\mu(x)\le \frac{\kappa}{2}\sum_{\e\in
E}|D_{\e}f|Q(\e),
\nneq
where ${\rm med}_{\mu}(f)$ is the median of $f$ under
$\mu$. The Cheeger inequality (\ref{func. Ch}) with the geometric
constant $\kappa$ is due to Diaconis and Stroock \cite{DS},  whose
idea goes back to Jerrum and Sinclair \cite{JS}. So Corollary
\ref{cor25} slightly improves  theirs in this particular case.

The Cheeger isoperimetric inequality for general jump  processes is studied by Chen and Wang \cite{CW}.
}\nrmk

\brmk{\rm If $G=(V,E)$ is a tree, i.e. there is only one path
without circle from $x$ to $y$ for any two different vertices
$x,y$, then the geometric constant $\kappa$ becomes {\it optimal}
for two types of metrics:

(a) $\rho(x,y)=1_{x\ne y}(\varphi(x)+\varphi(y))$ (then the constant
$\kappa$ in the weighted $L^1$-Poincar\'e inequality above is
optimal in the case of trees);

(b) $\rho(x,y)=\rho_w(x, y)$, the distance induced by some length function
$w$.

The optimality of $\kappa$ for the two types of metrics in the case
of trees is established by Liu-Ma-Wu \cite{LMW}, in a completely
different way. } \nrmk

The usual Cheeger inequality exhibits the relationship between
the isoperimetry and the  Poincar\'e inequality (\cite{AM,JS,LS}). Now we present the relationship between
the generalized Cheeger isoperimetric inequality and the Gaussian
concentration.

\bcor\label{cor38}   Assume  that $\sum_{y: y\sim x}q(x,y)\le B$ for all $x\in V$. Then
$$
c_{\rm G}\le \kappa^2 B.
$$

\ncor

\bprf This is due to \cite{GJLW}. But for the self-completeness, we still present
its proof. By the generalized Cheeger isoperimetric inequality in
Theorem \ref{thm-Cheeger} and the Cauchy-Schwarz inequality, we have
for any probability measure $\nu=f\mu$,
$$\aligned W_{1,\rho}(f\mu,\mu)
&\le \frac{\kappa}2\sum_{x\sim y}\mu(x)q(x,y)|f(x)-f(y)|\\
&\le \kappa \sqrt{I(\nu|\mu)}\sqrt{\frac12\sum_{x\sim
y}\mu(x)q(x,y)(\sqrt{f(x)}+\sqrt{f(y)})^2}\\
&\le  \kappa \sqrt{I(\nu|\mu)}\sqrt{2\sum_{x\sim
y}\mu(x)q(x,y) f(x)}\\
&\le  \kappa
\sqrt{I(\nu|\mu)}\sqrt{2B\sum_{x\in V}\mu(x)f(x)}=\kappa\sqrt{2BI(\nu|\mu)},
\endaligned$$
where the desired result follows. \nprf

\bcor\label{cor39} For the Laplacian $\LL=\Delta$ on the connected
graph $G=(V,E)$, we have for  the  graph metric $\rho_1$,
$$
 \kappa\le \frac{
d_*^2 b D}{|E|} ,
$$
where $d_*,D,b$ are given in Corollary \ref{cor24}.

 \ncor
\bprf Choosing $\gamma_{xy}$ distributed uniformly on all shortest
paths from $x$ to $y$ and $w=1$, since $Q(\e)=1/|E|, \mu(x)=d_x/|E|$, the geometric
constant $\kappa$ is bounded from above by
$$
\max_{\e\in E} |E|\ee^\gamma \sum_{x,y\in V} 1_{\gamma_{xy}}(\e)
\rho_1(x,y) \left(\frac{d_*}{|E|}\right)^2\le \frac{ d_*^2 b
D}{|E|}.
$$ \nprf

\section{Several  examples and graphs with symmetry}

\subsection{Several examples}
We begin with a baby-model.
 \bexa[Complete graph] {\rm Let $G=(V,E)$
be a complete graph with $n$ different vertices, i.e. for any different $x,y\in V$, $(x,y)\in E$
($n\ge 2$ of course). Consider the Laplacian $\LL=\Delta$ and the
graph metric $\rho_1$ which is now $\rho_1(x,y)=1_{x\ne y}$. Hence
$\mu$ is the uniform distribution on $V$ and $Q(\e)=1/|E|=
1/[n(n-1)]$. In such case the Dirichlet form is given by
$$
\EE(f,f) =\frac 12 \sum_{\e\in E} (D_ef)^2 Q(\e) = \frac {n}{n-1} {\rm
Var}_\mu(f),
$$
where ${\rm Var}_\mu(f)=\mu(f^2)-\mu(f)^2$ is the variance of $f$
w.r.t. $\mu$. So $c_{\rm P}=\frac {n-1}{n}$.

In this example, we take $\gamma_{xy}=\{(x,y)\}$ as random paths
and the length function $w\equiv 1$.

For the logarithmic Sobolev constant $c_{\rm LS}$, notice that for any fixed
edge $\e=(x_0,y_0)$, $\e\in \gamma_{xy}$ if and only if $x=x_0$ and
$y=y_0$; and $L_{w,\e}(x)=1_{x=x_0}/n$. Thus by Theorem \ref{thm11},
\begin{align}\label{LS complete}
 c_{\rm LS}&\le |E| \left({\rm Ent}_\mu(L_{w,\e})+\mu(L_{w,\e})\log(e^2+1)\right) \notag\\
&= n(n-1)\left[\frac 1{n^2}\log \frac 1n - \frac 1{n^2}\log \frac1{n^2} + \frac 1{n^2} \log(e^2+1) \right]\notag\\
&=(1-\frac 1n) \left[\log n + \log(e^2+1) \right]. \end{align}
Comparing with the  optimal logarithmic Sobolev constant $c_{\rm LS}=\frac{n-1}{n-2}\log (n-1)$ for complete graph (see \cite[Corollary A.5]{D-SC}),  the estimate \eqref{LS complete} has the correct order $\log n$.

\vskip0.3cm
Now we turn to bound the Gaussian concentration constant $c_{\rm G}$ by
the geometric quantity $K$ in (\ref{K}), w.r.t.   the graph metric
$\rho_1$. We have
$$
K= \max_{(x,y)\in E} \frac 1{Q(x,y)} \mu(x)\mu(y) =
\frac{|E|}{n^2}=\frac{n-1}{n}.
$$
Then by Theorem \ref{Thm WI}, for any $\mu$-probability density $f$,
$$
W^2_{1,\rho_1}(f\mu,\mu) =\frac 14 \left(\int_{V} |f-1| d\mu
\right)^2 \le 2 \frac{n-1}{n} \EE(\sqrt{f}, \sqrt{f})=2{\rm
Var}_\mu(\sqrt{f}).
$$
But by \cite[Theorem
3.1]{GLWY}, the corresponding optimal constant $c_{\rm G}=\frac{n-1}{2n}$. Consequently,  we have
$$
\frac K2=\frac{n-1}{2n}=c_{\rm G}\le K.
$$

\vskip0.3cm
For the generalized Cheeger isoperimetric inequality in Theorem
\ref{thm-Cheeger} w.r.t. $\rho_1$, we have
$$
c_{\rm I}\le \kappa=K\le \frac{n-1}{n}
$$
or equivalently
$$
\sum_{x\in V} |f(x)-\mu(f)|\mu(x) \le  \frac{n-1}{n} \sum_{\e\in
E}|D_{\e}f| Q(\e).
$$
This inequality becomes equality for indicator function $1_A$. Hence
$c_{\rm I}=\kappa$, i.e. our generalized Cheeger isoperimetric inequality
in Theorem \ref{thm-Cheeger} is optimal in this example.
 } \nexa

\begin{example}[Star] \label{star} {\rm Consider a star $G=(V,E)$ with a central vertex $v_0$ and $n$
outside vertices $\{v_i,i=1,\cdots,n\}$ connecting only with $v_0$.
For the Laplacian $\LL=\Delta$, we have
$\mu(v_0)=\frac{1}{2},\mu(v_k)=\frac{1}{2n}, 1\le k\le n$ and
$Q(\e)=\frac{1}{2n}$ for every edge $\e$. It is known that $c_{\rm P}=1$
(c.f.  \cite{DS}). Taking the length function $w\equiv 1$ in
(\ref{thm11b}), we obtain by calculus
$$
c_{\rm LS} \le \left(\frac32 -\frac 1n\right) \log [2n(e^2+1)].
$$
Applying the logarithmic Sobolev inequality to $f=1_{v_k}, 1\le k\le n,$ we get
$c_{\rm LS}\ge \log (2n)/2.$ Clearly, for large $n,$ we have the correct order $\log n.$

\vskip0.3cm
For the geometric constants $K$ and $\kappa$ associated with the
graph distance $\rho_1$, taking $\gamma_{xy}$ as the unique path
from $x$ to $y$ without circle, we have
$$
K\le \frac{9}{2}-\frac{4}{n} \,\,\text{and} \
\kappa=\frac{3}{2}-\frac{1}{n}.
$$
Considering $f(v)=2n{1}_{[v=v_1]}$,
 $$
W_1(f\mu,\mu)=\frac{3}{2}-\frac{1}{n}, \ \ \sum_{\e\in E} |D_{\e}f|Q(\e)=2,
 $$
we have that $$c_{\rm I}=\kappa=\frac{3}{2}-\frac{1}{n},$$ i.e., the geometric
quantity $\kappa$ as an upper bound of $c_{\rm I}$ is
 optimal.}
\end{example}

\bexa[Trees] {\rm Consider the full binary tree of depth $d$. For
$d\ge1$, such a tree has $2^{d+1}-1$ vertices, $2^{d+1}-2$ edges and
the maximum degree is $3$. Consider the Markov chain arising from
nearest neighbor random walk on this tree. The longest path is of
length $2d$ and  the value of $b$ defined in \eqref{b} is $(2^d-1)2^d$.  By Corollary
\ref{cor24}, we have $$K\le 18\cdot2^dd^3.$$
For $\kappa$ w.r.t. the graph distance $\rho_1,$ by \eqref{kappa}, Theorem 4.1 and Corollary 7.1 in \cite{LMW},
$$ \kappa=(2d-3)2^d+3\le  9\cdot d 2^{d-1},
$$ which shows that Corollary \ref{cor39} offers a good upper bound. } \nexa

\subsection{ Graphs with symmetry}
 In this section, we shall consider various graphs with symmetry.
 See   Chung \cite{C} for examples and  properties of symmetric
 graphs.
 For a graph $G=(V,E)$, an automorphism $f:V\rightarrow V$ is one-to-one mapping
 which preserves edges, i.e., for any $x,y\in V$, we have $(x,y)\in E$ if and only if $(f(x),f(y))\in E$.

 For any oriented edge $\e=(x,y)\in E$, consider the opposite oriented edge $\overleftarrow{\e}:=(y,x)$ and the non-oriented edge $\e^0:=\{x,y\}$. Put $E^0:=\{\e^0; \e\in E\}$, the set of all non-oriented edges.

\subsubsection{Edge-transitive graph}\label{sec edge} {\it A graph
$G$ is edge-transitive if, for any two  non-oriented  edges
$\{x,y\},\{x',y'\}$, there is an automorphism $f$ such that
$\{f(x),f(y)\}=\{x',y'\}$. }

\bcor\label{cor edge-tran} Assume that $G$ is edge-transitive. For
the Laplacian $\LL=\Delta$ and  the graph metric $\rho_1$, we have
$$
c_{\rm I}\le \kappa\le\ee \left[\rho_1^2(X,Y)\right]\ \ \text{and} \
\ c_{\rm G}\le \kappa^2 \le\left(\ee
\left[\rho_1^2(X,Y)\right]\right)^2,
$$
where the law of $(X,Y)$ is $\mu\times \mu$ and $\mu$ is the uniform measure on $V$.\ncor
 \bprf  We consider
a random (ordered) pair of vertices $(X,Y)$, chosen according to
$\mu\times\mu.$ Now given $(X,Y)$, we choose randomly a shortest
path $\gamma_{XY}$ between $X$ and $Y$ (uniformly chosen over all
possible shortest paths from $X$ to $Y$).

Notice that for $w\equiv 1$, $\rho=\rho_1$,
$$h(\e):=\frac{1}{Q(\e)w(\e)}\sum_{x,y\in V}\ee^\gamma
1_{\gamma_{xy}}(\e)\rho(x,y)\mu(x)\mu(y)=|E|\ee \rho_1(X,Y)
1_{\gamma_{XY}}(\e)$$ satisfies $h(\e)=h(\overleftarrow{\e})$ for $\e\in
E$ (that is true on any graph). Then $h$ can be regarded as a
function on $E^0$.

Now by the edge-transitivity, $h(\e)$ does not depend on $\e$, so we
get by averaging over all edges $\e$,
\begin{align*}
h(\e)&=\frac{1}{|E|} \sum_{\e\in E}|E|\ee
\left[\rho_1(X,Y)1_{\gamma_{XY}}(\e)
\right]\\
&=\ee \left[\rho_1^2(X,Y)\right].
\end{align*}Therefore
 $$\kappa\le \max_{\e\in E}h(\e)= \ee \left[\rho_1^2(X,Y)\right].
$$
Since $\sum_{y:y\sim x}q(x,y)\le1$ for all $x\in V$, we have $c_{\rm G}\le \kappa^2$ by Corollary \ref{cor38}. \nprf

\brmk\label{rmk correct order}{\rm  From the result above, we may
wonder whether on edge transitive graphs the correct order in diameter
$D$ of  $c_{\rm G}$ is $D^4$, and that
of $c_{\rm I}$ is  $D^2$, which is indeed true. For example, for the
Laplacian on the circle $\zz_p:=\zz/p\zz$, $c_{\rm P}$ is of order $D^2$. Taking the
eigenfunction $h$ corresponding to $\lambda_1=1/c_{\rm P}$ with
$\|h\|_{\rm Lip(\rho_1)}=1$, we see that ${\rm Var}_\mu(h)$ is of order
$D^2$, too. By the central limit theorem,
$$
\frac 1{\sqrt{t}} \int_0^t h(X_s) ds$$
converges weakly to the normal law $N(0,\sigma^2(h))$, where the
limit variance is
$$
\sigma^2(h)= 2\<(-\Delta)^{-1}h,h\>_\mu=2c_{\rm P} {\rm Var}_\mu(h),
$$
which is of order $D^4$. But from the Gaussian concentration
inequality we always have  $c_{\rm G}\ge \sigma^2(h)$. In other words $c_{\rm G}$
is at least of order $D^4$. We leave to the reader for verifying
that the correct order of $c_{\rm I}$ is $D^2$. }\nrmk

\bexa[Circle $\zz_p$]\label{circle} {\rm Let $p\in\zz$ and consider
the integers mod $p$ as $p$ points around a circle. For $x$ and $y$
in $\zz_p$, choose $\gamma_{xy}$ as the shorter of the two paths
from $x$ to $y$. For this model $\mu(x)=1/p \ (x\in \zz_p),
Q(\e)=1/(2p)$ for any edge $\e$.  It is well known that
$c_{\rm P}=(1-\cos\frac{2\pi}{p})^{-1}$, which is also the logarithmic Sobolev constant $c_{\rm LS}$ when $p$ is even (c.f. \cite{CC03}). Taking the length
function $w\equiv 1$ in (\ref{thm11b}),
we obtain by careful calculation,
\beq\label{circle upper} \aligned c_{\rm LS} &\le \begin{cases}\frac{\log(3(e^2+1))}{12} (p+1)(p+2), &\quad p\; \text{is even} ;\\\frac{\log(3(e^2+1))}{12} (p+1)(p+2)(1+3/p),  &\quad p
\; \text{is odd};\end{cases} \\
&\le 5 \log(3(e^2+1))(1-\cos\frac{2\pi}{p})^{-1},\endaligned \nneq
which, together with $c_{\rm LS}\ge c_{\rm P}=(1-\cos\frac{2\pi}{p})^{-1}$, offers a two sided estimate for $c_{\rm LS}$ with factor $5\log(3(e^2+1)).$
\vskip0.3cm
Since the circle $\zz_p$ is edge transitive, by Corollary \ref{cor
edge-tran}, we have
\begin{align*}
c_{\rm I}\le \kappa \le \frac{p^2}{12}+o(p^2);\ \ \ \ c_{\rm G}\le
\kappa^2\le\frac{p^4}{12^2}+o(p^4)
 .
\end{align*}
\vskip0.3cm
For this model,  Sammer and Tetali \cite{ST} proved that the best
constant in the transportation-entropy inequality \eqref{T1} is
$\frac{p^2}{48}+o(p^2)$. Then by Corollary \ref{cor23}, we have
$$
\frac{p^4}{24^2}+o(p^4)\le c_{\rm G}\le\frac{p^4}{12^2}+o(p^4).
$$
 Thus the correct order of $c_{\rm G}$ is $p^4$, as mentioned in Remark \ref{rmk correct order}.
}
\nexa

\subsubsection{Vertex-transitive graph}\label{sec vertex}{\it A graph $G$ is
vertex-transitive if, for any two vertices $u$ and $v$ there is an
automorphism $f$ such that  $f(u)=v$. } The automorphism group
defines an equivalent relation on the edges of $G$. Two undirected
edges $\e_1^0,\e_2^0$ are equivalent if and only if there is an
automorphism $\pi$ mapping $\e_1^0$ to $\e_2^0$. We can consider
equivalent classes of undirected edges, denoted by
$E^0_1,\cdots,E^0_s$.
  The index of $G$  is defined as
$$
{\rm index}(G)=\max_i \frac{|E^0|}{|E^0_i|}.
$$
 Clearly $|E^0_i|\ge |V|, i=1,\cdots,s,\ 1\le {\rm index}(G)\le d,$ where $d=d_x$ for any $ x\in V$, the degree of the graph
$G$.  See \cite{C}.

 For any edge $\e$ such that $\e^0\in E^0_i$, $i=1,\cdots, s$
\begin{align*}
\ee^{\mu\times\mu} \left[\rho_1(X,Y)1_{[\e\in\gamma_{XY}] }\right]
=&\frac{1}{2|E^0_i|} \sum_{\e: \e^0\in E^0_i}\ee^{\mu\times\mu} \left[\rho_1(X,Y)1_{[\e\in\gamma_{XY}] }\right]\\
\le &\frac{{\rm index}(G)}{|E|}\ee^{\mu\times\mu}
\left[\rho_1^2(X,Y)\right].
\end{align*}

For any two vertices $y,y'\in V$, there is an automorphism such that
$f(y)=y'$. Since the stationary distribution $\mu\equiv1/|E|$ and
$\rho_1(x, f(y))=\rho_1(f^{-1}(x),y)$ for any $x\in V$, we have
\begin{align*}
\ee^{\mu}\left[\rho_1^2(X,y')\right]=\ee^{\mu}\left[\rho_1^2(X,f(y))\right]
=\ee^{\mu}\left[\rho_1^2(f^{-1}(X),y)\right]=\ee^{\mu}\left[\rho_1^2(X,y)\right].
\end{align*}
Then for any fixed vertex $v_0\in V$,
$$\ee^{\mu\times\mu}\left[\rho_1^2(X,Y)\right]
=\ee^{\mu}\left[\rho_1^2(X,v_0)\right].$$ Therefore, we have proved
the following estimate of the Gaussian concentration constant
$c_{\rm G}$ in \eqref{WI-1}.
 \bcor\label{cor43} Let $\Delta$ be the Laplacian
operator on a vertex transitive graph $G$.  For any fixed vertex $v_0\in V$ and the
graph metric $\rho_1$, we have
 $$c_{\rm I}\le \kappa\le {\rm  index}(G)\
\ee^{\mu}\left[\rho_1^2(X,v_0)\right]\le
d\ee^{\mu}\left[\rho_1^2(X,v_0)\right]$$ and $$\ c_{\rm G}\le
K\wedge \kappa^2,$$ where $$ K\le{\rm index}(G)\
\ee^{\mu}\left[\rho_1^4(X,v_0)\right]\le
d\ee^{\mu}\left[\rho_1^4(X,v_0)\right] .$$

\ncor The proof of the bound on $K$ is similar, omitted here.

\subsubsection{Distance transitive graph} A
graph $G$ is distance transitive if, for any two pairs of vertices
$\{x,y\},\{x',y'\}$ with $\rho_1(x,y)=\rho_1(x',y')$, there is an automorphism
mapping $x$ to $x'$ and $y$ to $y'$. The distance transitive graph
is both edge-transitive and vertex-transitive. Then ${\rm index}(G)=1$. The estimate in
Corollary \ref{cor43} holds for
distance transitive graphs, that is, for any fixed vertex $v_0\in V$,
\beq\label{eq dis}c_{\rm G}\le  \left(\ee^{\mu}\left[\rho_1^2(X,v_0)\right]\right)^2\ \ \text{and } \ c_{\rm I}\le \ee^{\mu}\left[\rho_1^2(X,v_0)\right].\nneq

\begin{example}{\rm The vertices set $V$ consists of all the subsets of  $k$ elements in
$\{1,2,\cdots,n\}$ ($1\le k\le n$ is fixed). Define a metric on $V$
by $d(x,y)=k-|x\cap y|$. The edges set is given by $\{(x,y)\in V^2;
d(x,y)=1\}$. This is a distance transitive graph. The Markov process generated by the
Laplacian $\Delta$ on $G$ is known as the Bernoulli-Laplace
diffusion model. See Lee and Yau \cite{LY} for the estimate of the logarithmic Sobolev constant, Gao and Quastel \cite{GQ} for the exponential decay rate of entropy.

  By  \eqref{eq dis}, we have for the graph metric $\rho_1=d$,
$$c_{\rm I}\le \kappa\ \text{and}\ \ c_{\rm G}\le \kappa^2,$$
where
$$\kappa\le\frac{1}{\binom{n}{k}}\sum_{j=0}^{\min\{k, n-k\}} j^2
\binom{k}{j}\binom{n-k}{j}.
$$
}\end{example}

\end{document}